\documentclass[draft]{publmathd}
\usepackage{amsmath,amsfonts,amssymb}

\theoremstyle{plain}
\newtheorem{lemma}{Lemma}
\newtheorem{teo}{Theorem}
\newtheorem{corollary}{Corollary}

\DeclareMathOperator{\charac}{char}
\DeclareMathOperator{\expon}{exp}
\def\gp#1{\langle #1 \rangle}

\author{VICTOR BOVDI (Debrecen) and  ERNESTO SPINELLI (Lecce)}

\address{Victor Bovdi\\
Institute of Mathematics, University of Debrecen\\
H-4010 Debrecen, P.O.B. 12, Hungary\\
Institute of Mathematics and Informatics,\\ College of Ny\'\i
regyh\'aza\\ S\'ost\'oi \'ut 31/b, H-4410 Ny\'\i regyh\'aza,
Hungary} \email{vbovdi@math.klte.hu }

\address{Ernesto Spinelli\\
Dipartimento di Matematica "E. De Giorgi"\\
Universit\`{a} degli Studi di Lecce\\
Via Provinciale Lecce-Arnesano\\73100-LECCE\\
Italy}
\email{spinelli@ilenic.unile.it\\}

\title[Modular  group algebras with  maximal Lie nilpotency indices]
      {Modular  group algebras with \\ maximal Lie nilpotency indices}
\thanks{The research was supported by OTKA  No.T 037202 and  No.T 038059}
\dedicatory{Dedicated to the memory of Professor Jen\H o Erd\H os}



\submitted{.....}

\begin{document}

\begin{abstract}
In the present paper we give the full description of the Lie nilpotent
group algebras which  have maximal  Lie nilpotency indices.
\end{abstract}

\maketitle

\section{Introduction }

Let $R$ be an associative algebra with identity.  The algebra $R$
can be regarded as a Lie  algebra, called the associated Lie
algebra of $R$, via the  Lie  commutator $[x,y]=xy-yx$, for every
$x,y\in R$. Set
$$
[x_1,\ldots,x_{n-1},x_{n}]=[[x_1,\ldots,x_{n-1}],x_{n}],
$$ where
$x_1,\ldots,x_{n}\in R$. The \emph{$n$-th lower Lie power}
$R^{[n]}$ of $R$ is the associative ideal generated by all Lie
commutators $[x_1,\ldots,x_n]$, where  $R^{[1]}=R$ and
$x_1,\ldots,x_n \in R$. By induction, we define the \emph{$n$-th
upper  Lie power}  $R^{(n)}$ of $R$ as the associative ideal
generated by all Lie commutators $[x,y]$, where $R^{(1)}=R$ and
$x\in R^{(n-1)}$, $y\in R$.

An algebra  $R$ is called  \emph{Lie nilpotent} if there exists
$m$ such that $R^{[m]}=0$. The minimal integers $m,n$ such that
$R^{[m]}=0$ and $R^{(n)}=0$  are called \emph{ the  Lie nilpotency
index} and \emph{ the upper Lie nilpotency  index} of $R$ and
denoted they are by $t_{L}(R)$ and $t^{L}(R)$, respectively.

An algebra  $R$ is called  \emph{Lie hypercentral} if for any
sequence $\{a_i\}$ of elements of $R$ there exists some $n$, such
that $[a_1,\ldots,a_n]=0$.

Let $KG$ be the group algebra of a group $G$ over a  field $K$ of
characteristic $\charac(K)=p>0$. According to \cite{BKh,PPS} for
the noncommutative group algebra $KG$  the following statements
are equivalent: (a) $KG$ is Lie nilpotent (b) $KG$ is Lie
hypercentral; (c) $\charac(K)=p>0$, $G$ is nilpotent and its
commutator subgroup $G^{\prime}$ is a finite $p$-group. It is well
known \cite{P, SVB},  if $KG$ is Lie nilpotent then $t_{L}(KG)\leq
t^{L}(KG)\leq \vert G^{\prime}\vert +1$.  Moreover, according to
\cite{BP}, if $\charac(K)>3$, then $t_{L}(KG)= t^{L}(KG)$. But the
question when $t_{L}(KG)= t^{L}(KG)$ for $\charac(K)=2,3$ is still
open.

In the present paper we investigate  the group algebras $KG$ for
which  $t_{L}(KG)$ is {\it maximal}, i.e.    $t_{L}(KG)=\vert
G^{\prime}\vert +1$. In particular, if $G$ is a finite $p$-group
and $\charac(K)\geq 5$, then as Shalev \cite{S2} proved that
$\textsf{t}_{L}(KG)$ is maximal if and only if $G^{\prime}$ is
cyclic. We give a complete characterization by proving the
following:

\begin{teo}\label{main}
Let $KG$ be a Lie nilpotent group algebra  with
\newline $\charac(K)=p>0$. Then \quad $t_{L}(KG)=\vert G^{\prime}\vert
+1$ \quad if and only if one of the following conditions holds:
\begin{enumerate}
\item $G^\prime$ is cyclic;
\item $p=2$ and $G^\prime$ is the noncyclic of order $4$ and
$\gamma_{3}(G)\not=1$.
\end{enumerate}
\end{teo}

\begin{corollary}
Let $KG$ be a Lie nilpotent group algebra with
\newline $\charac(K)=p>0$. If \quad
 $t^{L}(KG)=\vert G^{\prime}\vert +1$, \quad then \quad
$t_{L}(KG)=t^{L}(KG)$.
\end{corollary}

By Du's Theorem \cite{D}, the previous result lists also the group
algebras $KG$ whose group of units $U(KG)$ has  maximal nilpotency
class under the assumption that $G$ is a finite $p$-group. Note
that Konovalov \cite{K} proved that $U(KG)$ of $2$-groups of
maximal class over field $K$ with  $\charac(K)=2$ have  the
maximal nilpotency class.

We use the  standard   notation for a group $G$:\quad $\Phi (G)$
denotes  the Frattini subgroup of $G$;\quad $g^h:=h^{-1}gh$ \quad
and  $(g,h):=g^{-1}h^{-1}gh$,\quad ($g,h\in G$); \quad
$\gamma_{i}(G)$ means  the $i$-{th} term of the lower central
series of $G$, i.e.
$$
\gamma_{1}(G)=G,\quad \quad
\gamma_{i+1}(G)=\big(\gamma_{i}(G),G\big)\quad \quad (i\geq 1).
$$
Moreover,   $C_n$ is the cyclic group of order $n$ and set
\begin{eqnarray*}
\mathrm{Q}_{2^n}&=&\gp{\; a,b\; \mid\;  a^{2^{n-1}}=1, \quad b^2=a^{2^{n-2}}, \quad  a^b=a^{-1}\; }, \quad\text{with}\quad n\geq 3;\\
\mathrm{D}_{2^n}& =&\gp{\; a,b\; \mid \; a^{2^{n-1}}=b^2=1,  \quad a^b=a^{-1}\; }, \quad\text{with}\quad n\geq 3;\\
\mathrm{SD}_{2^n}&=&\gp{\; a,b\; \mid\;  a^{2^{n-1}}=b^2=1, \quad  a^b=a^{{-1}+2^{n-2}}\; }, \quad\text{with}\quad n\geq 4;\\
\mathrm{MD}_{2^n}&=&\gp{\; a,b \; \mid\;  a^{2^{n-1}}=b^2=1, \quad a^b=a^{1+2^{n-2}}\; }, \quad\text{with}\quad n\geq 4.\\
\end{eqnarray*}

\section{Preliminaries}
Let $K$ be a field of characteristic $p>0$ and $G$ a group.
We consider  a sequence of subgroups of $G$ setting
$$
\mathfrak D_{(m)}(G)=G\cap ( 1+KG^{(m)}),\qquad\quad (m\geq 1).
$$
The subgroup $\mathfrak D_{(m)}(G)$ is  called the $m$-th   \emph{
Lie dimension subgroup} of $KG$. It is possible to describe  the
$\mathfrak D_{(m)}(G)$'s in terms of the lower central series of
$G$ in the following manner ( \cite{P}, p.44)
\begin{equation}\label{form}
\frak D_{(m+1)}(G)=\left \{
\begin{array}{ll}
G & \quad \text{if} \quad  m=0;\\
G^{\prime}& \quad \text{if}  \quad     m=1;\\
{\big(\frak D_{(m)}(G),G\big)(\frak D_{(\lceil {\frac{m}{p}}\rceil
+1)}(G))^p}& \quad \text{if}  \quad  m\geq 2.
\end{array} \right.
\end{equation}
where $\lceil {\frac{m}{p}}\rceil $ is the smallest integer greater than ${\frac{m}{p}}$.

Put $p^{d_{(m)}}:=[\frak D_{(m)}(G):\frak D_{(m+1)}(G)] $. If $KG$
is Lie nilpotent, according to Jennings' theory  \cite{S3} for the
Lie dimension subgroups, we get
$$
\displaystyle t^{L}(KG)=2+(p-1) \sum_{m\geq 1} md_{(m+1)}.
$$

\begin{lemma}\label{sha}(\cite{S1, S2})
Let $K$ be a field with $\charac(K)=p>0$ and $G$  a nilpotent
group such that  $G'$ is a finite $p$-group with
$\expon(G^{\prime})=p^l$.
\begin{enumerate}
\item  If $d_{(m+1)}=0$ and $m$ is a power of $p$, then
$\frak D_{(m+1)}(G)=1$.
\item If $d_{(m+1)}=0$ and $p^{l-1}$ divides
$m$, then $\frak D_{(m+1)}(G)=1$.
\end{enumerate}
\end{lemma}

\begin{lemma}\label{num}
Let $p, s, n\in \mathbb{N}$ and
 $m_0, \ldots, m_{s-1}$  the non-negative  integers  such that $s<n$ and
$\sum_{i=0}^{s-1}m_i =n$. Then \quad $\sum_{i=0}^{s-1}m_i p^i <
\sum_{i=0}^{n-1} p^i$.
\end{lemma}

\begin{proof} By the assumptions,  there exist integers $0\leq j_1 < \cdots <
j_k \leq s-1$ such that $m_{j_{l}}> 1$ for every $1\leq l\leq k$.
Since $p^s>p^{j_{k}}$, we obtain that
\begin{align*}
\sum_{i=0}^{s-1}m_{i}p^{i}=\sum_{i=0}^{s-1}p^{i}&+
\sum_{i=1}^{k}(m_{j_{i}}-1)p^{j_{i}}\\
&\leq \sum_{i=0}^{s-1}p^{i}+p^{j_{k}}(n-s)<
\sum_{i=0}^{s-1}p^{i}+\sum_{i=s}^{n-1}p^{i}.\quad
\end{align*}
\end{proof}

\begin{lemma}\label{L:3}
Let $K$ be a field with $\charac(K)=p>0$ and $G$ a nilpotent group
such that \quad $\vert G^{\prime} \vert=p^{n}$. Then \quad
$t^{L}(KG)=\vert G^{\prime} \vert +1$ if and only if \quad
$d_{(p^{i}+1)}=1$ and $d_{(j)}=0$, \quad where  $0\leq i\leq n-1$,
\quad $j\neq p^i +1$\quad and $j>1$.
\end{lemma}

\begin{proof} If $d_{(p^{i}+1)}=1$ for   $0\leq i\leq n-1$ and $d_{(j)}=0$
for $j>1$, then
\[
t^{L}(KG)=2+(p-1)\sum_{i=0}^{n-1}p^{i}=1+p^{n}=\vert
G^{\prime}\vert +1.
\]
In order to prove the other implication, we preliminarily remark
that
\begin{equation}\label{e:2}
\displaystyle \sum_{m\geq 2} d_{(m)}=n
\end{equation}
that is an immediate consequence of the definition of $d_{(j)}$'s.
Now we suppose that there exists $0\leq j\leq n-1$ such that
$d_{(p^{j}+1)}=0$. Let $s$ be the minimal integer for which
$d_{(p^{s}+1)}=0$. From  (\ref{form}) it follows at once that
$s\neq 0$ and by (1) of  Lemma \ref{sha} we have that $\frak
D_{(p^{s}+1)}(G)=1$ and so $d_{(r)}=0$ for every $r\geq p^{s}+1$.
It is immediate by (\ref{e:2}) that $\alpha :=
{\sum}_{i=0}^{s-1}d_{(p^{i}+1)} \leq n$. Let us consider the
following two cases: $\alpha =n$ and $\alpha < n$. If  $\alpha =n$
then, according to Lemma \ref{num}, we have that
\[
t^{L}(KG)=2+(p-1)\sum_{i=0}^{s-1}{p^{i}d_{(p^{i}+1)}}<
2+(p-1)\sum_{i=0}^{n-1}{p^{i}} =|G^\prime|+1.
\]
If $\alpha < n$  by (\ref{e:2}) there exists at least one
$j>1$ such that $d_{(j)}\neq 0$ and $j\neq p^{i}+1$. Suppose
$d_{(j_{1})},\ldots, d_{(j_{k})}$ are all of such  $d_{(j)}$'s,
where $j_{1}< \cdots < j_{k}$. Clearly,  $j_{k}\leq p^{s}$.
According to Lemma \ref{num} for the case $\alpha>s$, we obtain
\begin{align*}
t^{L}(KG)&=2+(p-1)\sum_{i=0}^{s-1}p^{i}d_{(p^{i}+1)}+
(p-1)\sum_{i=1}^{k}(j_{i}-1)d_{(j_{i})}  \\
&\leq 2+(p-1)\sum_{i=0}^{\alpha-1}{p^{i}}+
(p-1)(j_{k}-1)(n-\alpha)\\
&< 2+(p-1)\sum_{i=0}^{\alpha -1}{p^{i}}+(p-1)p^s(n-\alpha) \\
&< 2+(p-1)\sum_{i=0}^{\alpha -1}{p^{i}}+
(p-1)\sum_{i=\alpha}^{n-1}{p^{i}}=\vert G^{\prime} \vert +1.
\end{align*}
So, if $t^{L}(KG)$ is maximal, then  $d_{(p^{j}+1)}>0$ for each
$0\leq j\leq n-1$ and, by (\ref{e:2}), the Lemma is proved. \end{proof}

\begin{corollary} Let $K$ be a field with $\charac(K)=p>0$ and $G$
a nilpotent group with $\vert G^{\prime} \vert=p^{n}$. If \quad
$t^{L}(KG)=\vert G^{\prime} \vert +1$, then $\vert \frak
D_{(p^i+1)}(G)\vert =p^{n-i}$, \quad for $0\leq i\leq n$.
\end{corollary}

\begin{proof} By Lemma \ref{L:3}, it is easy to check that
\begin{align*}
\frak D_{(p^0+1)}(G)&\supset \frak D_{(p^1+1)}(G)\supset \frak
D_{(p^1+2)}(G)=\cdots\\
\cdots &=\frak D_{(p^i+1)}(G)\supset
\frak D_{(p^i+2)}(G)=\cdots \\
&=\frak D_{(p^s+1)}(G)\supset \frak D_{(p^s+2)}(G)=1
\end{align*}
for some $s\in \mathbb{N}$. Clearly,  $\vert \frak
D_{(p^s+1)}(G)\vert =p$, \quad $\vert \frak
D_{(p^{s-1}+1)}(G)\vert =p^2$ and $\vert \frak D_{(p^i+1)}(G)\vert
=p^{s-i+1}$, so  $\vert \frak D_{(p^0+1)}(G)\vert =p^{s+1}=p^n$
and $s=n-1$.\end{proof}

\begin{lemma}\label{L:4}
Let  $K$ be a field with  $\charac(K)=p>0$ and $G$  a nilpotent
group with $G^{\prime}$ a finite $p$-group such that
$t^{L}(KG)=\vert G^{\prime}\vert +1$.
\begin{enumerate}
\item If $p>2$ then $G^{\prime}$ is cyclic.
\item If $p=2$ then
$G^{\prime}$ has at most two generators.
\end{enumerate}
\end{lemma}

\begin{proof}
Assume that $\vert G^{\prime} \vert =p^n$. Let us prove that
\[
\vert \Phi (G^{\prime}) \vert \geq \left \{
\begin{array}{ll}
p^{n-1} &  if \quad  p\neq 2;\\
2^{n-2} &  if \quad   p=2.
\end{array} \right.
\]
First, set $p\neq 2$,\quad  $1<a<p$ and suppose that $\vert \Phi
(G^{\prime})\vert \leq  p^{n-2}$, where  $n\geq 2$. Since
$\expon(G^{\prime}/\Phi (G^{\prime}))=p$, we have that
$\expon(G^{\prime})=p^k\leq p^{n-1}$ for some $k$. By Lemma
\ref{L:3}, we get $d_{(ap^{n-2}+1)}=0$ and $p^{k-1}$ divides
$p^{n-2}$. Then by (2) of Lemma \ref{sha} we obtain that $\frak
D_{(ap^{n-2}+1)}(G)=1$. But $ap^{n-2}<p^{n-1}$ and  $\frak
D_{(p^{n-1}+1)}(G)\neq 1$, which is a contradiction.

Now, set $p=2$ and we suppose that $\vert \Phi(G^{\prime})\vert
\leq 2^{n-3}$, where  $n\geq 3$. By   Lemma \ref{L:3} we have
$d_{(3\cdot 2^{n-3}+1)}=0$. Since $3\cdot 2^{n-3}< 2^{n-1}$, by
(2) of Lemma \ref{sha} we get that $\frak D_{(3\cdot
2^{n-3}+1)}(G)=1$ and $\frak D_{(2^{n-1}+1)}(G)\neq 1$, which is a
contradiction either. \end{proof}

\begin{lemma}\label{L:5}
Let $K$ be a field with  $\charac(K)=2$,  $G$ a nilpotent group
such that $G^\prime$ is  $2$-generated   finite $2$-group
and let \quad $t^{L}(KG)=\vert G^{\prime}\vert +1$. If either
$\gamma_2(G)^2\subset \gamma_3(G)$ or
$\gamma_3(G)\cap\gamma_2(G)^2=1$ then $\vert \gamma_3(G)\vert =2$ and
$\gamma_2(G)\cong C_2\times C_2$.
\end{lemma}

\begin{proof}
Assume that $\vert G^{\prime} \vert =2^n$.
Let $G$ be nilpotent of class $cl(G)=t\leq n+1$ and
$\gamma_2(G)^2\subset \gamma_3(G)$. Then  by Theorem III.2.13
(\cite{H}, p.266), we have that $\gamma_k(G)^2 \subseteq
\gamma_{k+1}(G)$ for every $k\geq 2$. Let us prove by induction on
$i$ that $\frak D_{(2^i+1)}(G)= \gamma_{i+2}(G)$.
It follows at once that $\frak D_{(2)}(G)=\gamma_{2}(G)$ and
$\frak D_{(3)}(G)=\gamma_{3}(G)$. According to Lemma
\ref{L:3} we have
\begin{eqnarray*}
\frak D_{(2^{i+1}+1)}(G)&=&\frak D_{(2^i+2)}(G)=\big(\frak
D_{(2^i+1)}(G),G\big)\cdot\frak D_{(\lceil
2^{i-1}+\frac{1}{2}\rceil+1)}(G)^2\\
&=&\big(\gamma_{i+2}(G),G\big)\cdot\frak D_{(
2^{i-1}+2)}(G)^2=\gamma_{i+3}(G)\cdot\frak D_{( 2^{i-1}+2)}(G)^2\\
&=&\gamma_{i+3}(G)\cdot\frak D_{(
2^{i}+1)}(G)^2=\gamma_{i+3}(G)\cdot\gamma_{i+2}(G)^2\\
&=&\gamma_{i+3}(G).
\end{eqnarray*}
It follows  that $\frak D_{(2^{t-1}+1)}(G)=\gamma_{t+1}(G)=1$, but
by Lemma \ref{L:3} we have  $\frak D_{(2^{n-1}+1)}(G)\not=1$,\quad
so $t>n$ and  $t=n+1$.

 Obviously, for \quad $i\geq 1$
\begin{eqnarray*}
\frak D_{( 2^{i}+2)}(G)&=&\big(\frak D_{(
2^{i}+1)}(G),G\big)\cdot\frak D_{(
2^{i-1}+2)}(G)^2=\gamma_{i+3}(G)\cdot\frak D_{( 2^{i}+1)}(G)^2\\
&=&\gamma_{i+3}(G)\cdot\gamma_{i+2}(G)^2\\
&=&\gamma_{i+3}(G);\\
\frak D_{( 2^{i}+3)}(G)&=&\big(\frak D_{(
2^{i}+2)}(G),G\big)\cdot\frak D_{(
2^{i-1}+2)}(G)^2\\
&=&\gamma_{i+4}(G)\cdot\gamma_{i+2}(G)^2.\\
\end{eqnarray*}
Since $\frak D_{( 2^{i}+2)}(G)=\frak D_{( 2^{i}+3)}(G)$ for $i\geq
1$ we get that
$$
\gamma_{i+3}(G)=\gamma_{i+4}(G)\cdot\gamma_{i+2}(G)^2.
$$
According to   $\gamma_{3}(G)^2\supset\gamma_{4}(G)^2\supset\cdots$
it  follows that
$$
\gamma_{4}(G)=\gamma_{3}(G)^2\cdot\gamma_{4}(G)^2\cdot
\gamma_{5}(G)^2\cdots\gamma_{t}(G)^2=\gamma_{3}(G)^2.
$$
Since $\Phi(\gamma_{3}(G))=\gamma_{3}(G)^2$ we have
$$
[\gamma_{3}(G):\Phi(\gamma_{3}(G))]=[\gamma_{3}(G):\gamma_{4}(G)]=2,
$$
so $\gamma_{3}(G)$ is cyclic. According to Theorem 12.5.1 in
\cite{Ha} the $2$-generated group $\gamma_{2}(G)$ with cyclic
subgroup of index $2$ is one of the following groups:
$\mathrm{Q}_{2^n}$, $\mathrm{D}_{2^n}$, $\mathrm{SD}_{2^n}$,
$\mathrm{MD}_{2^n}$, or $C_2\times C_{2^{n-1}}$, and therefore
$\gamma_{2}(G)^2=\gamma_{4}(G)$.

 Moreover,
$\gamma_{3}(G)^2\subseteq \gamma_{5}(G)$. Indeed, the elements of
the form $(x,y)$, where $x\in \gamma_2(G)$ and $y\in G$ are
generators of $\gamma_3(G)$, so we have to  prove that $(x,y)^2
\in \gamma_5(G)$. Evidently,
$$
(x^2,y)=(x,y)(x,y,x)(x,y)=(x,y)^2(x,y,x)^{(x,y)}
$$
and  $(x^2,y),\; (x,y,x)^{(x,y)}\in \gamma_5(G)$, so $(x,y)^2\in
\gamma_5(G)$ and $\gamma_3(G)^2\subseteq\gamma_5(G)$.
Using the  fact that
$exp(\gamma_{3}(G)\big/\gamma_{5}(G))=2$, since $\gamma_{3}(G)$ is
cyclic, we obtain that $\vert \gamma_{3}(G)\vert =2$ and $\gamma_2(G)\cong
C_2\times C_2$.

Now,  let $\gamma_3(G)\cap\gamma_2(G)^2=1$. By (\ref{form}) we have that
$\frak{D}_{(2)}(G)=\gamma_2(G)$,
$\frak{D}_{(3)}(G)=\gamma_2(G)^2\cdot\gamma_3(G)$ and
\[
\begin{split}
\frak{D}_{(2)}(G)/\frak{D}_{(3)}(G)&=\gamma_2(G)/\big[\gamma_2(G)^2\cdot\gamma_3(G)\big]\\
&\cong \big[\gamma_2(G)/\gamma_2(G)^2\big] \big/\gamma_3(G).
\end{split}
\]
Since \quad $|\frak{D}_{(2)}(G)/\frak{D}_{(3)}(G)|=2$ and
$|\gamma_2(G)/\gamma_2(G)^2|=4$, from the last equality it follows
that $|\gamma_3(G)|=2$ and $\gamma_4(G)=1$.

 Obviously,
$(\gamma_2(G),\gamma_2(G))\subseteq \gamma_4(G)=1$, so
$\gamma_2(G)$ is abelian and
\begin{eqnarray*}
\frak{D}_{(4)}(G)=\frak{D}_{(5)}(G)&=&\big(\gamma_2(G)^2\cdot\gamma_3(G),G\big)(\gamma_2(G)^2\cdot\gamma_3(G))^2\\
&=&\gamma_2(G)^4\cdot\gamma_3(G)^2=\gamma_2(G)^4.\\
\end{eqnarray*}
It is easy to check that
\[
\Phi(\gamma_2(G)^2\cdot\gamma_3(G))=(\gamma_2(G)^2\cdot\gamma_3(G))^2=\gamma_2(G)^4=\frak{D}_{(5)}(G).
\]
Therefore $\gamma_2(G)^2\cdot\gamma_3(G)$ is a cyclic subgroup of
index $2$ in $\gamma_2(G)$ and
\[
\gamma_2(G)=\gp{a}\times\gp{b}\cong C_{2^{n-1}} \times
C_2\quad\quad(|a|=2^{n-1},\; |b|=2).
\]
Clearly, $\gamma_2(G)^2=\gp{a^2}$ and either $\gamma_3(G)=\gp{b}$
or $\gamma_3(G)=\gp{a^{2^{n-2}}b}$. Now, let us compute the weak
complement of $\gamma_3(G)$ in $G^{\prime}$ (see \cite{BK}, p.34).
It is easy to see that $\nu(b)=\nu(a^{2^{n-2}}b)=2$ and the weak
complement will be $A=\gp{a}$. Since $G$ is of class $3$, by (ii)
of Theorem 3.3  (\cite{BK}, p.43) we have
\begin{eqnarray*}
t_L(KG)=t^L(KG)=2^{n}+1&=&
t(\gamma_2(G))+t(\gamma_2(G)/A)\\
&=&2^{n-1}+3,
\end{eqnarray*}
so $n=2$ and $\gamma_2(G)\cong C_{2} \times C_2$.\end{proof}

\section{ Proof of   Theorem \ref{main}}

Let $KG$ be a Lie nilpotent group algebra  with $\charac(K)=p>0$
and let \quad $t^{L}(KG)=\vert G^{\prime}\vert +1$. By Lemma
\ref{L:4} is either $p>2$ and $\gamma_2(G)$ is cyclic or $p=2$ and
$\gamma_2(G)$ has at most $2$ generators.

Now, let $p=2$ and  $\gamma_2(G)$  a   $2$-generated group. Let us
prove that either $\gamma_2(G)^2\subset\gamma_3(G)$ or
$\gamma_3(G)\cap\gamma_2(G)^2=1$.

First, we suppose that $\gamma_3(G)\subseteq \gamma_2(G)^2$. It is easy
to see, that
$$
\frak{D}_{(2)}(G)=\gamma_2(G),\quad\quad
\frak{D}_{(3)}(G)=\gamma_2(G)^2
$$
and
\[
\begin{split}
\frak{D}_{(2)}(G)/\frak{D}_{(3)}(G)&=\gamma_2(G)/\gamma_2(G)^2\\
&\cong \gamma_2(G)/\Phi(\gamma_2(G))\cong C_2\\
\end{split}
\]
which contradicts to the fact that  $\gamma_2(G)$  is a $2$-generated group.

Finally, suppose  $\gamma_3(G)\cap\gamma_2(G)^2\not=1$ and
$\gamma_2(G)^2\not\subset\gamma_3(G)$. Clearly,
\begin{eqnarray*}
\frak{D}_{(2)}(G)&=&\gamma_2(G);\\\qquad\qquad
\frak{D}_{(3)}(G)&=&\gamma_2(G)^2\cdot \gamma_3(G);\\
\frak{D}_{(2^i+1)}(G)&\equiv &\gamma_3(G)^{2^{i-1}}\cdot
\gamma_2(G)^{2^{i}}\pmod{\gamma_4(G)}\quad\quad (i\geq 2).
\end{eqnarray*}
Since $\frak{D}_{(4)}(G)\equiv
\frak{D}_{(5)}(G)\pmod{\gamma_4(G)}$, it follows that
\begin{eqnarray*}
2&=&[\frak{D}_{(3)}(G):\frak{D}_{(4)}(G)]\\
&\equiv&
[\gamma_2(G)^2\cdot
\gamma_3(G):(\gamma_2(G)^2\cdot \gamma_3(G))^2]\\
&\equiv& [(\gamma_2(G)^2\cdot \gamma_3(G)):\Phi(\gamma_2(G)^2\cdot
\gamma_3(G))]\pmod{\gamma_4(G)}\\
\end{eqnarray*}
and  $L\equiv \gamma_2(G)^2\cdot \gamma_3(G)\pmod{\gamma_4 (G)}$
is a cyclic group. Set $L=\gp{a}$. Since
$[\gamma_2(G):\gamma_2(G)^2]=4$, we get that  $L$ is a  subgroup
of index $2$ in $\gamma_2(G)/\gamma_4(G)$. Thus, by Theorem 12.5.1
of \cite{Ha} we have that
$$
\gamma_2(G)/\gamma_4(G)=\gp{\;a,b\;\mid\;a^{2^{n-1}}=1,\;
b^2\in\gp{a}\;\;}
$$
is one of the following groups:\quad $C_{2^{n}}$,
$C_{2^{n-1}}\times C_{2}$, $\mathrm{Q}_{2^n}$, $\mathrm{D}_{2^n}$,
$\mathrm{SD}_{2^n}$ or $\mathrm{MD}_{2^n}$. But in these cases
$\gamma_2(G)^2/\gamma_4(G)= \gp{a^2}$ and from
$$
L\equiv\gamma_2(G)^2\cdot\gamma_3(G)\pmod{\gamma_4(G)}
$$
it follows  that $\gp{a}\equiv
\gp{a^2}\cdot\gamma_3(G)\pmod{\gamma_4(G)}$, so
$a\in\gamma_3(G)/\gamma_4(G)$ and $ \gamma_2(G)^2\subset
\gamma_3(G)\pmod{\gamma_4(G)}$,  a contradiction. Therefore, by
Lemma \ref{L:5}, we have that $|\gamma_3(G)|=2$ and $\gamma_2(G)\cong
C_2\times C_2$.

Conversely, let $p\geq 2$ and $\gamma_2(G)$  a cyclic group. By
(ii) of Theorem 3.1 of \cite{BK} we get $t_{L}(KG)=t^{L}(KG)=\vert
G^{\prime}\vert +1$. Now, let $p=2$ and $G^\prime$  the noncyclic
group  of order $4$ and $\gamma_{3}(G)\not=1$. Again, by (ii) of
Theorem 3.3 of \cite{BK} we obtain that
$$
t_{L}(KG)=t^{L}(KG)=t(G^\prime)+t(G^\prime/A)=\vert
G^{\prime}\vert +1,
$$
where $A$ is the weak complement of $\gamma_3 (G)$ in $G^\prime$
and the proof is complete.

\end{document}